# A NEW COEXISTENCE RESULT FOR COMPETING CONTACT PROCESSES

BY BENJAMIN CHAN AND RICHARD DURRETT

*Cornell University*

Neuhauser [*Probab. Theory Related Fields* **91** (1992) 467–506] considered the two-type contact process and showed that on $\mathbb{Z}^2$ coexistence is not possible if the death rates are equal and the particles use the same dispersal neighborhood. Here, we show that it is possible for a species with a long-, but finite, range dispersal kernel to coexist with a superior competitor with nearest-neighbor dispersal in a model that includes deaths of blocks due to "forest fires."

**1. Introduction.** In this paper we consider a two-type contact process in which the state of each site $x \in \mathbb{Z}^2$ is $0 =$ vacant, or 1 and 2 that indicate occupancy by two different species of trees. We formulate the evolution as follows:

(i) Particles of type $i$ die at rate $\delta_i$, give birth at rate $\beta_i$.

(ii) A 1 born at $x$ is sent to a $y$ chosen according to a truncated power-law distribution $p_1(x, y) = c_1 \mathbb{1}_{\{\|y-x\|_\infty = 1\}} + c_2 \mathbb{1}_{\{1 < \|y-x\|_\infty \leq M\}} \|y - x\|_\infty^{-\rho}$, where $c_1, c_2 > 0$ and $\rho < 3$.

(iii) A 2 born at $x$ is sent to a $y$ chosen at random from $\mathcal{N}_1 = \{y : \|x - y\|_\infty = 1\}$, that is, for these $y$, $p_2(x, y) = 1/8$.

(iv) If $y$ is already occupied, then the birth is suppressed.

(v) For each $x$, death to $\{y : \|x - y\|_\infty \leq F/2\}$ occurs at rate $\delta_0$.

The rules above imply that each type of particle behaves individually like a contact process, but in addition there are deaths to blocks of size $F \times F$ due to forest fires.

Neuhauser [7] studied the model with $\delta_0 = 0$ and in which both species have the same dispersal distribution. She proved that coexistence is impossible on $\mathbb{Z}^2$ if both species have the same death rate and conjectured that in general the superior competitor (the species with the higher reproductive ratio $\beta_i/\delta_i$) would win the competition. In this paper, we show that it is possible for two species to coexist in a model with forest fires if the weaker competitor has larger dispersal range.

This result is, in a sense, an application of ecology to mathematics since it is based on a phenomenon that is well known to occur in forests; see [8]. The strongest competitors, which ultimately take over an undisturbed patch of forest, are those that need the least amount of light to grow. However, other species that are better at colonizing gaps in the forest created by disturbances can coexist with them.

Our new result is not a counterexample to Neuhauser's conjecture, nor should it cast doubt on it. The competitive exclusion principle, which has been proved in the setting of ordinary differential equation models by Levin [5], states that the number of coexisting species cannot exceed the number of resources they compete for. (To be precise, Levin's result is for generic parameter values. In Neuhauser's model if $\beta_1 = \beta_2$ and $\delta_1 = \delta_2$, coexistence is possible in $d \geq 3$, but in this case we do not really have two species.) In Neuhauser's situation the two species compete for space, so there should only be one winner, even if different dispersal distributions are allowed. In our setting newly disturbed space is a second resource, so a species adapted to that niche can coexist with a superior competitor.

To produce our example we begin by fixing $\beta_2$, $\delta_2$ and the ratio $\beta_1/\delta_1$ so that

$$(1) \qquad \beta_2/(1 + \beta_1/\delta_1) > \delta_2 \lambda_c(\mathcal{N}_1), \qquad \beta_1/\delta_1 > \lambda_c(\mathcal{N}_1),$$

where $\lambda_c(\mathcal{N}_1)$ is the critical value of $\beta/\delta$ for the contact process with neighborhood $\mathcal{N}_1$.

THEOREM 1. *Suppose* (1). *We can choose* $\delta_1$, $\delta_0$, $F$ *and* $M$ *so that* 1's *and* 2's *coexist.*

The second condition in (1) implies that when $c_2 = 0$ particles of type 1 can survive in the absence of 2's. To explain the first condition, note that 1's die at rate $\delta_1$ and are born at a site at rate $\leq \beta_1$ so the set of sites occupied by 1's is dominated by a product measure with density $\beta_1/(\delta_1 + \beta_1)$. In the limit as $\delta_1 \to \infty$ with $\beta_1/\delta_1$ fixed, the environments seen by the 2's at different birth attempts are independent in the dominating process. Thus in the limit as $\delta_1 \to \infty$ the 2's dominate a contact process with birth rate $\beta_2 \delta_1/(\beta_1 + \delta_1)$ and death rate $\delta_2$. Thus it suffices to have $\beta_2/(1 + \beta_1/\delta_1) > \delta_2 \lambda_c(\mathcal{N}_1)$ for the 2's to survive in the system with no fire. Once this is done it is easy to use known results for the contact process to show if $\delta_0 \leq c/F^2$ and $c$ is small, 2's also survive in the presence of fires.

PROPOSITION 1. *Suppose* (1). *There are constants $c$ and $\Delta$ so that for any $F > 0$, the 2's survive whenever $\delta_1 \geq \Delta$ and $\delta_0 \leq c/F^2$.*

The 1's survive by migrating from existing patches to newly created ones. By choosing the fire size $F$ appropriately, the probability that a square contains a cleared area which remains open (free of 2's) for a long period of time is high. We choose $\delta_0$ to produce enough such gaps, but also to not burn them again before we are ready.

PROPOSITION 2. *Suppose* (1). *Fix $\delta_1 \geq \Delta$. If $\delta_0 = a/F^3$ with $a$ sufficiently small, $F \geq F_0(a)$ and $M$ is chosen appropriately, then the 1's survive.*

Propositions 1 and 2 are proved by block constructions, so together they imply the existence of a stationary distribution in which both types are present with positive density. The remainder of the paper is devoted to proving the two propositions.

**2. Preliminaries on contact processes.** In this paper, unless stated otherwise, we denote a contact process on $\mathbb{Z}^2$ with neighborhood $\mathcal{N}_1$ and initial configuration $A$ by $A_t^A$. We let $_LA_t^A$ be the corresponding process truncated to not allow births outside $[-L, L]^d$, and suppose that $A \subset [-L, L]^d$. Define $\tau^A = \min\{t : A_t^A = \varnothing\}$ to be the extinction time for $A_t^A$, and $_L\tau^A$ the extinction time for $_LA_t^A$. We now state two results for the contact process that will be used later on. First a survival condition, due to Bezuidenhout and Grimmett [1], which we take from Theorems 2.12 and 2.23 of [6].

THEOREM 2. *$A_t$ survives if and only if it satisfies the following condition:*

*For every $\varepsilon > 0$ there are choices of $n$, $L$, $T$ so that*

$$P\{_{L+2n}A_{T+1}^{[-n,n]^d} \supset x + [-n,n]^d \text{ for some } x \in [0,L)^d\} > 1 - \varepsilon$$

*and*

$$P\{_{L+2n}A_{t+1}^{[-n,n]^d} \supset x + [-n,n]^d \text{ for some } 0 \leq t < T,$$
$$\text{and for some } x \in \{L+n\} \times [0,L)^{d-1}\} > 1 - \varepsilon.$$

Note that in contrast to the usual block construction, one of the new copies of the completely filled square is on the side of the space–time box.

The next theorem is a shape theorem proved by Durrett and Griffeath [3]. For the contact process they did this only for large enough birth rate. Given the results described in Chapter 2 of [6], it is routine to extend this to all supercritical contact processes. Readers who do not want to take this leap of faith can replace $\lambda_c(\mathcal{N}_1)$ by the critical value for the one-dimensional nearest-neighbor contact process.

Let $\eta_t$ be a supercritical contact process, which includes the Richardson model (contact process with death rate 0) as a special case. Let $\nu$ be the upper invariant measure of $\eta_t$, which is the limit starting from all 1's, and let $\eta_t^\nu$ be the stationary process starting from $\nu$. Let $\tau = \min\{t : \eta_t^0 = \varnothing\}$ and $t(x) = \min\{t : x \in \eta_t^0\}$. The set of sites hit by time $t$ is

$$H_t = \{y \in \mathbb{R}^d : \exists x \in \mathbb{Z}^d \text{ with } \|x - y\|_\infty \le 1/2 \text{ and } t(x) \le t\}.$$

The region where $\eta_t^0$ and $\eta_t^\nu$ are coupled is

$$K_t = \{y \in \mathbb{R}^d : \exists x \in \mathbb{Z}^d \text{ with } \|x - y\|_\infty \le 1/2 \text{ and } \eta_t^0(x) = \eta_t^\nu(x)\}.$$

THEOREM 3. *Let $\eta_t^0$ be either a supercritical contact process or a Richardson process (contact process with death rate 0) starting from 0. There is a convex set $\mathcal{S} \subset \mathbb{R}^d$ such that for any $\varepsilon > 0$,*

$$(1 - \varepsilon)\mathcal{S} \subset t^{-1} H_t \subset (1 + \varepsilon)\mathcal{S}$$

*and*

$$(1 - \varepsilon)\mathcal{S} \subset t^{-1}(H_t \cap K_t) \subset (1 + \varepsilon)\mathcal{S}$$

*for $t \ge T_0(\varepsilon)$ a.s. on $\{\tau = \infty\}$, where $T_0(\varepsilon) < \infty$ is a random time.*

**3. Proof of Proposition 1.** We begin by constructing our process from a collection of Poisson processes. For $x, y \in \mathbb{Z}^2$ and $i = 1, 2$, let $\{T_n^{i,x,y} : n \ge 1\}$, $\{U_n^{i,x} : n \ge 1\}$ and $\{V_n^x : n \ge 1\}$ be the arrival times of Poisson processes with rates $\beta_i p_i(x, y)$, $\delta_i$ and $\delta_0$. At times $U_n^{i,x}$ we put a $\delta_i$ at $x$ to kill the particle at $x$ if it is of type $i$. At times $T_n^{i,x,y}$ we draw an arrow of type $i$ from $x$ to $y$ to indicate that there is a birth from $x$ to $y$ if $x$ is occupied by type $i$ and $y$ is empty. At time $V_n^x$ we put a $\delta$ at all $y$ such that $\|y - x\|_\infty \le F/2$ to kill all the particles inside the square. Our process has finite range so it follows from a result of Harris [4] that this gives a process well defined for all time.

We first show that the 2's survive when there are no forest fires, that is, $\delta_0 = 0$. This is done by showing that the set of sites occupied by 2's can outcompete a process that is stochastically larger than the sites occupied by 1's. Let $\tilde{\xi}_t^i$ be the set of sites occupied by $i$'s when $\delta_0 = 0$. Let $\tilde{\zeta}_t^i$ be the set of sites occupied by $i$'s when $\delta_0 = 0$ and each $x \in \mathbb{Z}^2$ not in state 2 flips from 0 to 1 at rate $\beta_1$ and 1 to 0 at rate $\delta_1$ independent of all other sites.

LEMMA 1. *Suppose $\tilde{\xi}_0 = \tilde{\zeta}_0$. Using $\le$ for stochastically smaller than, we have $\tilde{\xi}_t^1 \le \tilde{\zeta}_t^1$.*

PROOF. Since the death rate of 1's in $\tilde{\xi}_t$ is equal to that of $\tilde{\zeta}_t$ while the birth rate of 1's is smaller, this follows from Theorem 1.5 in Chapter 3 of [6]. □

LEMMA 2. *Suppose $\beta_2$, $\delta_2$ and the ratio $r = \beta_1/\delta_1$ are given. Let $\gamma = \beta_2\delta_1/(\beta_1+\delta_1)$ and let $A = \tilde{\zeta}_0^2$. Writing a second superscript to indicate the value of $\beta_2$, given $\varepsilon, \theta > 0$ there is a $\Delta$ such that there is a coupling with*

$$P({}_L A_t^{A,\gamma-\varepsilon} \leq {}_L\tilde{\zeta}_t^2 \leq {}_L A_t^{A,\gamma+\varepsilon} \text{ for } 0 \leq t \leq T) > 1 - \theta$$

*whenever $\delta_1 \geq \Delta$ and $A \subset \mathbb{Z}^2 \cap [-L,L]^2$.*

PROOF. In the graphical representation, each site $x$ receives a type-2 arrow at rate $\beta_2$. The arrivals in $[0,T]$ that touch some $x \in [-L,L]^2$ are a Poisson point process. Let $T_0 = 0 < T_1 < T_2 < \cdots$ be the arrival times and let $N(t)$ be the number of arrivals by time $t$. If $N$ is large, $P(N(t) \geq N) \leq \theta/3$. Having chosen this value of $N$, if $\alpha$ is small,

$$P\left(\min_{1 \leq i \leq N} T_i - T_{i-1} < \alpha\right) \leq \theta/3.$$

Since $\eta_t$ is a Markov process, the effect of a type-2 arrow to $y$ at time $t+h$ depends only on the state of the process at time $t$. The process that flips between 1 and 0 at rates $\delta_1$ and $\beta_1$ is a two-state Markov process, so if $\delta_1$ is large and $\beta_1/\delta_1$ is fixed, the total variation distance between the process starting from 1 or 0 and product measure is smaller than $\theta/3N$ at all times $t \geq \alpha$. Since the process of 2-arrows is independent of the 1 flipping, and the latter has a stationary distribution $\beta_1/(\beta_1+\delta_1)$, combining the last three estimates proves the result. □

Note that this implies $\tilde{\zeta}_t^2 \Rightarrow A_t^{A,\gamma}$ for all finite sets $A$ and $t < \infty$, where $\Rightarrow$ denotes weak convergence. Also by construction one has $\tilde{\xi}_t^2 \geq \tilde{\zeta}_t^2$ for all $t \geq 0$. Therefore once we show that $\tilde{\zeta}_t^A$ survives the forest fires, the proof is complete. But this is easy if we limit the frequency of our fires so that they are but small perturbations on the graphical representation.

PROOF OF PROPOSITION 1. For any chosen $\varepsilon$, the idea is to choose $\delta_0$ small enough so that the chance of a forest fire interfering with a given space–time block has probability $\leq \varepsilon$ and hence can be neglected.

Let $A_t$ be a contact process with birth rate $\gamma' = \beta_2\delta_1/(\beta_1+\delta_1) - \varepsilon'$ and death rate $\delta_2$. Since $\beta_2\delta_1/(\beta_1+\delta_1) > \delta_2\lambda_c(\mathcal{N}_1)$, we can pick $\varepsilon'$ small enough to have $\gamma' > \delta_2\lambda_c(\mathcal{N}_1)$. Then by Theorem 2, there are choices of $n$, $L$ and $T$ so that in the absence of fires

$$P\{{}_{L+2n}A_{T+1}^{[-n,n]^2} \supset x + [-n,n]^2 \text{ for some } x \in [0,L]^2\} > 1 - \varepsilon$$

and

$$P\{_{L+2n}A_{t+1}^{[-n,n]^2} \supset x + [-n,n]^2 \text{ for some } 0 \le t < T,$$
$$\text{and for some } x \in \{L+n\} \times [0,L)\} > 1 - \varepsilon.$$

Since a fire eliminates all the points within $[x - F/2, x + F/2]^2$, our space–time block $[-L - 2n, L + 2n]^2 \times [0, T+1]$ is affected only if the center of a fire is inside $[-F/2 - L - 2n, F/2 + L + 2n]^2$. Let $\tau = (F + 2L + 4n + 1)^2(T+1)$. The probability our space–time block is unaffected by fires is bounded below by

$$1 - \int_0^\tau \delta_0 e^{-\delta_0 t}\,dt = e^{-\delta_0 \tau} > 1 - \varepsilon$$

if $\delta_0 \le c/F^2$, where $c > 0$ is some constant. Lemma 2 states that $_{L+2n}\tilde{\zeta}_t^{2,[-n,n]^2} \ge {}_{L+2n}A_t^{[-n,n]^2}$ with probability $> 1 - \varepsilon$ whenever $\delta_1 \ge \Delta$. Combining the estimates above we get

$$P\{_{L+2n}\zeta_{T+1}^{2,[-n,n]^2} \supset x + [-n,n]^2 \text{ for some } x \in [0,L)^2\} > (1-\varepsilon)^3$$

and

$$P\{_{L+2n}\zeta_{t+1}^{2,[-n,n]^2} \supset x + [-n,n]^2 \text{ for some } 0 \le t < T,$$
$$\text{and for some } x \in \{L+n\} \times [0,L)\} > (1-\varepsilon)^3.$$

Setting $\varepsilon$ small enough, we have recreated the block events of Theorem 1. For reasons that are explained in detail in Liggett's [6] book, $\zeta_t^2$ dominates a supercritical $k$-dependent oriented percolation and therefore has a positive probability of survival. This completes the proof since by construction $\xi_t^2 \ge \zeta_t^2$ for all time. □

## 4. Proof of Proposition 2.

4.1. *Size of a fire box.* We first determine an appropriate size $F$ for our fires. Recall that the 1's survive by colonizing the gaps created by forest fires. So we want a fire to create enough space for the 1's to migrate from afar, and survive within the area for a long period of time. To accommodate these needs, we define three squares centered at 0, labeled $B_1$, $B_2$ and $B_3$ with sides of lengths $2L$, $2r_1 L$ and $2r_2 L = 2F$, where $1 < r_1 < r_2$. The innermost square $B_1$ provides the space for 1's to land, the square $B_2$ allows them room to grow and the perimeter area $B_3 - B_2$ serves as a buffer against the 2's. We want to show the following.

PROPOSITION 3. *Suppose $\beta_2$, $\delta_2$ are fixed and $c_1 \beta_1 > \delta_1 \lambda_c(\mathcal{N}_1)$. Here $c_1$ is the coefficient for short-range births in the dispersal function $p_1(x,y)$,* and $\lambda_c(\mathcal{N}_1)$ is the critical value of $\beta/\delta$ for survival of the contact process with neighborhood $\mathcal{N}_1$. Then there are choices of $r_1$, $r_2$ and $T$ so that in the absence of interference of other forest fires the following events happen with probability $\geq 1 - \theta$ for large $L$:

(i) If $B_1$ receives a type-1 particle whose descendants survive within $B_2$ for a period of $T$, then they will survive within $B_2$ for a period of $7T/2$.

(ii) If $B_3$ does not contain any 2's at time 0, 2's do not reach $B_2$ by time $7T/2$.

In the first lemma we prove (i). To grow the set of 1's we only use the short-range edges, and hence have a version of the process $A_t^x$. Recall that $r_1 L \tau$ is the extinction time of $r_1 L A_t^x$. Let $H_t^x$ and $K_t^x$, as in Theorem 3, be the set of sites hit and coupled for $A_t^x$. For simplicity let $C_t^x = H_t^x \cap K_t^x \subset H_t^x$. Let $D_R(y)$ be a square of side $2R$ centered at $y$ and let $D_R = D_R(0)$. Since $\beta$ and $\delta$ are fixed, there are $R_{1,2} > R_{1,1} > 0$ such that $D_{R_{1,1}} \subset \mathcal{S} \subset D_{R_{1,2}}$, where $\mathcal{S}$ is the limiting shape of $A_t^0$. Let

$$(2) \quad G_t^x = \{(1-\varepsilon)D_{R_{1,1}} \subset t^{-1}(C_t^x - x), t^{-1}(H_t^x - x) \subset (1+\varepsilon)D_{R_{1,2}}\}.$$

LEMMA 3. *Let $A_t$ be a supercritical contact process with birth rate $\beta$ and death rate $\delta$. Let $T = 8L/7R_{1,1}$ and let $x \in [-L, L]^2$. For any $\theta > 0$, we can choose $r_1$ so that for large $L$*

$$P(A_t^x \subset B_2 \text{ for } t \leq 7T/2 |_{r_1 L}\tau^x \geq T) \geq 1 - \theta$$

*and*

$$P(_{r_1 L}\tau^x \geq 7T/2 |_{r_1 L}\tau^x \geq T) \geq 1 - \theta.$$

PROOF. Let $\varepsilon = 1/2$ in (2). Theorem 2 implies that there is a $t^* > 0$ such that $P(G_t^x \text{ for all } t \geq t^* | \tau^x = \infty) \geq 1 - \theta$. For $L \geq t^* R_{1,1}/4$, we have on $E = \bigcap_{t \geq t^*} G_t^x$,

$$C_{4L/R_{1,1}}^x \supset x + [-2L, 2L]^2 \supset [-L, L]^2,$$

$$H_{4L/R_{1,1}}^x \subset tx + [-K, K]^2 \quad \text{where } K = \tfrac{3}{2}R_{1,2}(4L/R_{1,1}).$$

Letting $r_1 = 6R_{1,2}/R_{1,1} + 1$ we have on $E$,

$$P([-L, L]^2 \subset C_{4L/R_{1,1}}^x \subset H_{4L/R_{1,1}}^x \subset D_{r_1 L} | \tau^x = \infty) > 1 - \theta.$$

Now $\{G_{7T/2}^x, \tau^x = \infty\} \subset \{G_{7T/2}^x, \tau^x \geq T\}$, so

$$P(G_{7T/2}^x | \tau^x \geq T) \geq P(G_{7T/2}^x | \tau^x = \infty) \frac{P(\tau^x = \infty)}{P(\tau^x \geq T)}.$$

The last fraction tends to 1 as $T \to \infty$. By the definitions of $G_t^x$ and $r_1$,
$$\{\tau^x \geq T\} \cap G_{7T/2}^x = \{r_1 L \tau^x \geq T\} \cap G_{7T/2}^x.$$
Since $P(r_1 L \tau^x \geq T) \leq P(\tau^x \geq T)$ we have (recall $7T/2 = 4L/R_{1,1}$)
$$P(G_{7T/2}^x | \tau^x \geq T) \leq P(G_{7T/2}^x | r_1 L \tau^x \geq T)$$
$$\leq P(r_1 L \tau^x \geq 7T/2 | r_1 L \tau^x \geq T)$$
and the desired conclusions follow. $\square$

To prepare for the proof of Lemma 5 we observe that $P(G_{7T/2}^x | \tau^x = \infty) \to 1$ as $T \to \infty$ so working backward through the last few calculations shows that for large $L$
$$(3) \qquad P(r_1 L \tau^x \geq 7T/2) \geq (1 - \theta) P(\tau^x = \infty).$$

The next lemma proves (ii) in Proposition 3.

LEMMA 4. *If $r_2$ is chosen large enough, then, for large $L$, $A_0 \cap B_3 = \varnothing$ implies*
$$P(A_t \cap B_2 = \varnothing \text{ for all } t \leq 7T/2) \geq 1 - \theta.$$

PROOF. Fix $\varepsilon = 1/2$. Let $R_t^A$ be the Richardson process with neighborhood $\mathcal{N}_1$ and growth rate $\beta_2$. Let $\mathcal{S}$ be the convex set that is the limit in Theorem 3. Suppose $R_{2,2} > R_{2,1} > 0$ are two numbers such that $D_{R_{2,1}} \subset \mathcal{S} \subset D_{R_{2,2}}$. Let
$$G_t = \{(1 - \varepsilon) D_{R_{2,1}} \subset t^{-1} R_t^0 \subset (1 + \varepsilon) D_{R_{2,2}}\}.$$
There is a constant $t^*$ such that $P(G_t \text{ for all } t \geq t^*) \geq 1 - \delta$.

Suppose $t_2 > t_1 \geq t^*$. Then
$$P(G_{t_1}) \geq P(G_{t_1}, t_2^{-1} R_{t_2}^0 \subset (1 + \varepsilon) D_{R_{2,2}}) \geq 1 - \delta.$$
Since $R_t$ is a Markov process,
$$\delta > P(G_{t_1}, t_2^{-1} R_{t_2}^0 \not\subset (1 + \varepsilon) D_{R_{2,2}})$$
$$\geq P(t_2^{-1} R_{t_2 - t_1}^{(1-\varepsilon) t_1 D_{R_{2,1}}} \not\subset (1 + \varepsilon) D_{R_{2,2}}) P(G_{t_1}).$$
Since $P(G_{t_1}) \geq 1 - \delta$, taking $t_2 = t_1 + t$ gives
$$P((t + t_1)^{-1} R_t^{(1-\varepsilon) t_1 D_{R_{2,1}}} \not\subset (1 + \varepsilon) D_{R_{2,2}}) \leq \delta/(1 - \delta).$$
Let $\theta = \delta/(1 - \delta)$. Let $m \geq r_1 L$ and $m \geq (1 - \varepsilon) t^* R_{2,1}$. We can find $t_1 \geq t^*$ so that $m = (1 - \varepsilon) t_1 R_{2,1}$ and as a result
$$P(R_{7T/2}^{[-m,m]^2} \subset (3/2)(7T/2 + t_1) D_{R_{2,2}}) \geq 1 - \theta.$$

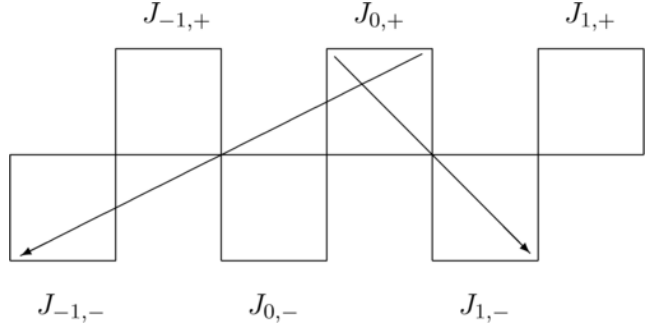

Fig. 1. *Boxes in the block construction.*

Let $r_2 = (3/2)(7T/2 + t_1)R_{2,2}/L$ and $U = D_{r_2 L}$. Since $R_t$ is self-dual and $B_2 \subset [-m,m]^2$,

$$P(R_{7T/2}^{U^c} \cap B_2 = \varnothing) = P(R_{7T/2}^{B_2} \cap U^c = \varnothing)$$
$$= P(R_{7T/2}^{B_2} \subset U)$$
$$\geq 1 - \theta. \qquad \square$$

4.2. *Block construction.* Like Proposition 1, Proposition 2 will also be shown by using a block construction. We will make use of the following notation. Let $\mathcal{L} = \{(m,n) \in \mathbb{Z}^2 : m+n \text{ is even}\}$ and $e_1 = (1,0)$. We have already defined $F = r_2 L$. Let $\alpha > 0$, $k = F^\alpha$, $W = F + F^{1+\alpha}$, and define

$$J_- = [-kW, 0]^2, \qquad J_+ = (0, kW]^2, \qquad J_{m,\pm} = 2mkWe_1 + J_\pm.$$

We also define

$$B(i,j) = [iW, (i+1)W] \times [jW, (j+1)W], \qquad i,j \in \mathbb{Z},$$
$$B_2(x,y) = [x - r_1 L, x + r_1 L] \times [y - r_1 L, y + r_1 L].$$

We say there is a *source of 1's* at $B(i,j)$ at stage $n$ if there is a block $B_2(x,y) \subset B(i,j)$ that is empty of 2's and contains at least one type-1 particle at all times between $2nT$ and $2(n+1)T$, where $T$ will be chosen later. We say there is a *source of 1's* in $J_{m,+}$ at stage $n$ if there is a source of 1's in one of the $B(i,j)$'s that make up $J_{m,+}$.

The next lemma is the key to proving Proposition 2.

LEMMA 5. *Suppose we have a source of 1's at $J_{0,+}$ at stage 0. $k$ and $W$ can be chosen so that with probability $\geq 1 - \theta$, there are sources of 1's at $J_{-1,-}$ and $J_{1,-}$ at stage 1.*

We alternate between $+$ and $-$ for our sources to avoid possible interference. Once Lemma 5 is proven, we will choose the cutoff $M = 4kW$. The events described in Lemma 5 will have a finite range of dependence, so a standard argument (see, e.g., [2]) implies that our process will dominate a supercritical oriented percolation on $\mathcal{L}$, and Proposition 2 follows.

PROOF OF LEMMA 5. There are three steps in creating a source in an adjacent box:

(a) fire clears an area in $[0, T/2]$ and no fire touches the clearing in $[T/2, 4T]$;

(b) dispersal brings a 1 to the box $B_1$ inside the clearing during $[T/2, T]$;

(c) the immigrant survives to time $4T$ without leaving $B_2$, and without 2's coming into $B_2$.

We will show that the intersection of these events has probability $\geq \eta > 0$ in a given $W \times W$ square in $J_{-1,-}$ and $J_{1,-}$. Since the $J$'s are a grid of $k^2$ such squares, if $k$ is large, success will have high probability. At this point the reader might worry about $\eta$ becoming small as $k$ gets large but this point will be addressed in the proof.

The probability of the first desired event happening in one specified $W \times W$ area is at least

$$\int_0^{(W-F)^2 T/2} \delta_0 e^{-\delta_0 t}\, dt \cdot \int_{(2F)^2 7T/2}^{\infty} \delta_0 e^{-\delta_0 t}\, dt$$
$$= (1 - e^{-\delta_0 (W-F)^2 T/2}) e^{-14 \delta_0 F^2 T}.$$

If we choose $W = F + F^{1+\alpha}$, $\alpha > 0$, then to get

$$e^{-14\delta_0 F^2 T} \geq 1 - \theta \quad \text{and} \quad 1 - e^{-\delta_0(W-F)^2 T/2} \geq 1 - \theta$$

we need

$$\delta_0 F^2 T \leq -\log(1-\theta)/14 \quad \text{and} \quad \delta_0 F^{2+2\alpha} T/2 \geq -\log(\theta).$$

The first is guaranteed by our assumption $\delta_0 \leq c/F^3$ for a suitable choice of $c$. Once we fix $\delta_0$ the second will hold when $L$ and hence $F$ is large.

Suppose there is a source in $J_{0,+}$. It contains at least one particle for all times between $T/2$ and $T$. Then the maximum $L^\infty$ distance between our source and an area cleared by a nice fire in $J_{1,-}$ or $J_{-1,-}$ (see Figure 1) is $4kW$. Thus, the minimum rate at which our source spreads to the target $B_1$ within the cleared area is $u = c_2 \beta_1 |4kW|^{-\rho} L^2$. If we let $k = F^\alpha$, then there is a $c < \infty$ so that the probability of spread to the target area between $T/2$ and $T$ is bounded below by

$$1 - \exp(-uT/2) \geq 1 - \exp(-c\beta_1 L^{3-(1+2\alpha)\rho}).$$

To check this note that $W \leq 2F^{1+\alpha}$ and recall $F$ and $T$ are multiples of $L$. If $(1+2\alpha)\rho < 3$, then for large $L$ the last quantity is $\geq 1 - \theta$.

Once a particle has landed on an empty area, the probability that its descendant will survive there until $t = 4T$ is bounded below by

$$P(_{r_1L}A^x_{7T/2} \neq \varnothing) \geq (1-\theta)P(\tau^x = \infty)$$

as shown in (3).

Combining the last three estimates, we see that the probability of producing a source in a given $W \times W$ area in $J_{1,-}$ of $J_{-1,-}$ is bounded below as $L \to \infty$. The number of independent opportunities, $k = F^\alpha$, tends to $\infty$ so the probability of success tends to 1 and the proof is complete. $\square$

DEPARTMENT OF MATHEMATICS
CORNELL UNIVERSITY
105 MALOTT HALL
ITHACA, NEW YORK 14853
USA
E-MAIL: chanb@math.cornell.edu

DEPARTMENT OF MATHEMATICS
CORNELL UNIVERSITY
523 MALOTT HALL
ITHACA, NEW YORK 14853
USA
E-MAIL: rtd1@cornell.edu
URL: www.math.cornell.edu/~durrett